\documentclass[10pt]{article}
\usepackage{amssymb,amsmath,pstricks}
\usepackage{graphics,graphicx}
\usepackage{setspace}

\makeatletter\@addtoreset {equation}{section}\makeatother

\def\bexa{\begin{example}}\def\eexa{\eex\end{example}}
\def\brem{\begin{remark}}\def\erem{\eex\end{remark}}
\def\bthm{\begin{theorem}}\def\ethm{\end{theorem}}
\def\blem{\begin{lemma}}\def\elem{\end{lemma}}
\def\bcor{\begin{corollary}}\def\ecor{\end{corollary}}
\def\bdefi{\begin{definition}}\def\edefi{\end{definition}}
\def\beq{\begin{equation}}\def\eeq{\end{equation}}

\newcommand{\R}{{\mathbb R}}
\newcommand{\N}{{\mathbb N}}

\def\CN{{\cal N}}

\def\noi{\noindent}
\def\pa{{\partial}}

\newcommand{\bi}{\begin{itemize}}\newcommand{\ei}{\end{itemize}}
\newcommand{\ben}{\begin{enumerate}}\newcommand{\een}{\end{enumerate}}
\newcommand{\bce}{\begin{center}}\newcommand{\ece}{\end{center}}

\def\eps{\varepsilon}

\newcommand{\barr}{\begin{array}}\newcommand{\earr}{\end{array}}
\newcommand{\bpm}{\begin{pmatrix}}\newcommand{\epm}{\end{pmatrix}}
\newcommand{\ba}{\begin{array}}\newcommand{\ea}{\end{array}}
\def\ri{{\rm i}}

\def\del{\delta}

\def\ba{\begin{array}} \def\ea{\end{array}}
\def\bgat{\begin{gather}} \def\egat{\end{gather}}
\def\eps{\varepsilon}

\def\eex{\hfill\mbox{$\rfloor$}}

\def\al{\alpha}

\def\paxpml{\partial_x^\text{{\tiny PML}}}
\def\paypml{\partial_y^\text{{\tiny PML}}}
\def\pazpml{\partial_z^\text{{\tiny PML}}}
\def\upml{u^\text{{\tiny PML}}}
	\def\uref{u^\text{{\tiny REF}}}
\def\uhatpml{\hat{u}^\text{\tiny{PML}}}

\newtheorem{theorem}{Theorem}[section]
\newtheorem{definition}[theorem]{Definition}

\newtheorem{lemma}[theorem]{Lemma}

\newtheorem{remark}[theorem]{Remark}
\newtheorem{corollary}[theorem]{Corollary}
\newtheorem{example}[theorem]{Example}

\begin{document}

\title{\bf Perfectly Matched Layers for Coupled Nonlinear Schr\"{o}dinger Equations with Mixed Derivatives}

\author{Tom\'{a}\v{s} Dohnal\footnote{Email: dohnal@math.uka.de, Tel: +497216087670, Fax: +497216086679}
\\
  {\small Institut f\"{u}r Angewandte und Numerische Mathematik 2,
    Universit\"{a}t Karlsruhe, Kaiserstr. 12, 76128 Karlsruhe, Germany}}

\date{\today}
\maketitle

\begin{abstract}
  This paper constructs perfectly matched layers (PML) for a system of 2D  Coupled Nonlinear Schr\"{o}dinger equations with mixed derivatives which arises in the modeling of gap solitons in nonlinear periodic structures with a non-separable linear part. The PML construction is performed in Laplace Fourier space via a modal analysis and can be viewed as a complex change of variables. The mixed derivatives cause the presence of waves with opposite phase and group velocities, which has previously been shown to cause instability of layer equations in certain types of hyperbolic problems. Nevertheless, here the PML is stable if the absorption function $\sigma$ lies below a specified threshold. The PML construction and analysis are carried out for the linear part of the system. Numerical tests are then performed in both the linear and nonlinear regimes checking convergence of the error with respect to the layer width and showing that the PML performs well even in many nonlinear simulations.
  
  \end{abstract}

\noi  


\section{Introduction}
Perfectly matched layers (PML) are a relatively simple and efficient tool for the truncation of spatial domains for wave type problems posed on unbounded (or large) domains. In numerical simulations such a truncation is often necessary as well as desired and PML guarantees that waves traveling through the boundary are absorbed and reflections that are only exponentially small with respect to the layer width occur. PML have been first proposed by B\'{e}renger \cite{Beren94} for Maxwell's equations and since then derived and analyzed for many other equations, like wave and Helmholtz equations \cite{TY98}, linearized Euler equations \cite{H98,H03}, general first order hyperbolic systems \cite{AHK06}, Schr\"{o}dinger equation \cite{C97,H03_b}, etc. This paper proposes, analyzes and tests PML for a 2D Coupled Nonlinear Schr\"{o}dinger system (CNLS) with mixed derivatives
\beq\label{E:CME}
\ri\pa_t u_j +(\al^{(x)}_j\pa_x^2+\al^{(y)}_j\pa_y^2+\beta_j\pa_x\pa_y)u_j +\Gamma \CN_j(u_1,\ldots,u_N)=0, \quad j \in \{1,\ldots, N\},
\eeq
where $\CN_j$ is a polynomial (typically cubic or cubic-quintic) nonlinearity and $\al^{(x)}_j,\al^{(y)}_j, \beta_j$ and $\Gamma$ are real numbers. In \cite{DU09}, where it is called a system of Coupled Mode Equations, this system is shown to be an asymptotic model for gap solitons in the 2D periodic nonlinear Schr\"{o}dinger equation with a finite contrast non-separable periodic potential. Previously the author together with T. Hagstrom have studied in \cite{DH07} PML for 1D and 2D Coupled Mode Equations governing gap solitons in periodic structures with infinitesimal contrast. In that work the modal analysis in Laplace-Fourier space was used. The 1D problem was hyperbolic and the general PML construction for first order hyperbolic systems \cite{AHK06} employing auxiliary variables was used. The 2D case was of a mixed type and required a combination of the method in \cite{AHK06} and a complex coordinate stretching. The problem \eqref{E:CME} at hand is of generalized Schr\"{o}dinger type and the modal analysis in Laplace-Fourier space reveals that a complex change of coordinates is sufficient for PML construction. In this paper \eqref{E:CME} is studied under the condition $\al_j^{(x)}\al_j^{(y)}>\beta_j^2  \quad \forall j \in \{1,\ldots, N\}$, which is dictated by the asymptotic derivation in \cite{DU09} and implies ellipticity of the spatial linear operator in  \eqref{E:CME}.

As  \eqref{E:CME} is a model for nonlinear solitary waves, the following scenarios are particularly relevant for numerical investigations: evolution of a perturbed solitary wave; interaction of a solitary wave with a defect in the medium; or collision of several solitary waves. All these processes will typically lead to shedding of radiation that usually has small amplitude compared to the pulse(s), travels away from them, and needs to be treated at the boundary of the domain of interest. Examples of studies using PML in such situations are \cite{BEVK05,DA05,GW08}. In addition, simulations where a pulse of magnitude comparable to the solution maximum leaves the domain are often desired. In such a case the polynomial nonlinearity cannot be in general neglected in the layers and leads to truly nonlinear layer dynamics.

The derivation and analysis of PML in this paper is based purely on the \textit{linear} part of the system \eqref{E:CME}, nevertheless the presented numerical tests demonstrate satisfactory functionality even in prototypical examples corresponding to all of the above \textit{nonlinear} scenarios. The analysis guarantees that in the linear regime the layer is absorbing and perfectly matched. Stability of the layer equations in time is shown to hold if the maximum of the absorption function $\sigma$ lies below a threshold, which diverges to infinity for $\beta \rightarrow 0$. The layer equations are, therefore, conditionally stable, which is in spite of the presence of plane waves in the linear part of \eqref{E:CME} with \textit{opposite group and phase velocities}.  Such a mismatch of group and phase velocities has been shown in certain hyperbolic systems to lead to instability of the layer equations \cite{BFJ03,AK06}. PML for 3D linear Schr\"{o}dinger equations with mixed derivatives have been previously used in \cite{Cheng_etal07}. Perfect matching and stability were, however, not analyzed there in the presence of mixed derivatives.

The rest of the paper is organized as follows. In Section \ref{S:CNLS_relevance} relevance of the CNLS system is discussed and conditions on coefficients are provided which allow removal of the mixed derivatives via a change of variables. The PML is derived via the modal analysis in Laplace-Fourier space in Section \ref{S:derive_PML}. In Section \ref{S:stab} stability of the linear ($\Gamma=0$) layer equations is analyzed. Finally, Section \ref{S:num_tests} presents a number of numerical tests in both the linear and nonlinear regimes. Exponentially fast convergence of the error within the physical domain with respect to the layer width is verified via linear tests and observed, though with larger error values, even in nonlinear tests.

\section{Relevance of the CNLS with Mixed Derivatives}\label{S:CNLS_relevance}
Systems of the type \eqref{E:CME} have been shown in \cite{DU09} to describe gap solitons in ``nonseparable'' Kerr nonlinear structures for values of their spectral parameter (frequency or propagation constant) lying in an asymptotic neighborhood of a spectral gap edge. The CNLS system is then called the Coupled Mode Equations (CMEs) and governs the dynamics of slowly varying envelopes of the gap soliton. The particular model for which CMEs were derived in \cite{DU09} was the periodic nonlinear Schr\"{o}dinger equation 
\beq\label{E:PNLS}
\ri \pa_{t'} \psi+\Delta \psi - V(x',y')\psi +\Gamma |\psi|^2 \psi=0, \quad V(x'+d_1,y')=V(x',y'+d_2)=V(x',y') \ \forall (x',y')\in \R^2
\eeq
with some $d_{1,2}>0$ and where the periodic structure $V$ is fixed, i.e., does not depend on the asymptotic parameter, and non-separable. In physics literature this case is typically referred to as a large contrast periodic structure. Equation \eqref{E:PNLS} describes propagation of light in 2D photonic crystals as well as evolution of matter waves in Bose Einstein condensates. The nonseparability condition requires that $V(x',y') \neq V_1(x')+V_2(y')$ for any functions $V_1, V_2$. For a gap soliton near a gap edge defined by $N$ maxima or minima of the corresponding band structure $\omega(k)$ the general CMEs read as \eqref{E:CME} with the coefficients $\al^{(x)}_j, \al^{(y)}_j$ and $\beta_j$ proportional to the second derivatives of the spectral bands at the extrema with respect to the components of the wavevector $k$, see \cite{DU09}. In particular, $\beta_j$ is the mixed second derivative. As the extrema are either all minima or all maxima, the coefficients satisfy $\text{sign}(\al^{(x)}_j)=\text{sign}(\al^{(y)}_j) \ \forall j\in \{1, \ldots, N\}$ and 
\beq\label{E:coef_cond}
\al^{(x)}_j\al^{(y)}_j>\beta_j^2 \quad \forall j\in \{1, \ldots, N\},
\eeq
which also guarantees ellipticity of the spatial differential operator in \eqref{E:CME}.

In contrast, structures with a separable linear part, like the periodic nonlinear Schr\"{o}dinger equation \eqref{E:PNLS} with $V(x,y) = V_1(x)+V_2(y)$, which  was studied in \cite{SY07,DPS09}, lead to $\beta_j=0 \ \forall j$ and the CMEs take the form of classical CNLS systems.

CMEs of the type \eqref{E:CME} can also be derived as an approximative model in the same asymptotic regime as above for gap solitons in the Maxwell problem $\Delta \psi - V(x',y') \partial_{t'}^2 \psi - \Gamma \partial_{t'}^2 (\psi^3)=0, \ \psi(x,t)\in \R$, with a finite contrast periodic structure $V(x',y')$ as these do not have a separable linear part either.

In \cite{DU09} an example of the potential $V(x',y')$ is presented, for which the band structure indeed leads to $\beta_j\neq 0$.

A simple prototypical example of the system \eqref{E:CME}, which is used in this paper for some of the numerical tests, is
\beq\label{E:CME2}
\begin{split}
\ri\partial_t u_1 + (\al^{(x)}_1\pa_x^2+\al^{(y)}_1\pa_y^2+\beta_1\pa_x\pa_y)u_1 + \Gamma \left[|u_1|^2u_1 +(2|u_2|^2u_1+u_2^2\bar{u}_1)+\eps_q|u_1|^4u_1\right]=&0\\
\ri\partial_t u_2 + (\al^{(x)}_2\pa_x^2+\al^{(y)}_2\pa_y^2+\beta_2\pa_x\pa_y)u_2 + \Gamma \left[|u_2|^2u_2 +(2|u_1|^2u_2+u_1^2\bar{u}_2)+\eps_q|u_2|^4u_2\right]=&0,
\end{split}
\eeq 
where the quintic nonlinearity with $\eps_q<0$ has been included in order to avoid issues with blowup of solutions of the cubic 2D nonlinear Schr\"{o}dinger equation \cite{SS00}.

\paragraph{Removal of the mixed derivative via a change of variables}

The mixed derivatives in \eqref{E:CME} (and in \eqref{E:CME2}) can be in certain cases removed via the change of variables
\beq\label{E:transf}
\bpm\tilde{x}\\\tilde{y}\epm = \bpm a\cos \theta & -b\sin \theta \\ a \sin \theta & b \cos \theta\epm \bpm x \\ y\epm
\eeq
for some $a,b\in\R$, leading to
\[
\begin{split}
\al^{(x)}_j\pa_x^2+\al^{(y)}_j\pa_y^2+\beta_j\pa_x\pa_y = &\left(\al^{(x)}_ja^2\cos^2\theta  +\al^{(y)}_jb^2\sin^2\theta - \beta_j {ab\over 2}\sin 2\theta \right)\pa_{\tilde{x}}^2 \\
&+\left(\al^{(x)}_ja^2\sin^2\theta  +\al^{(y)}_jb^2\cos^2\theta + \beta_j {ab\over 2}\sin 2\theta \right)\pa_{\tilde{y}}^2 \\
&+\left((\al^{(x)}_ja^2-\al^{(y)}_jb^2)\sin 2\theta +\beta_j a b \cos 2\theta\right)\pa_{\tilde{x}}\pa_{\tilde{y}}.
\end{split}
\]
The mixed derivative is, clearly, removed if $\tfrac{\al^{(x)}_1}{\al^{(y)}_1}=\tfrac{\al^{(x)}_2}{\al^{(y)}_2}= \ldots = \tfrac{\al^{(x)}_N}{\al^{(y)}_N}$ by the choice $a=b \left({\al^{(x)}_j \over \al^{(y)}_j}\right)^{1/2}$ and $\theta=\frac{\pi}{4}$. Otherwise, the removal is successful if there are constants $a,b \in \R$ such that
\beq\label{E:remove_cond}
\tfrac{\beta_1}{\al^{(x)}_1a^2-\al^{(y)}_1b^2}=\tfrac{\beta_2}{\al^{(x)}_2a^2-\al^{(y)}_2b^2} = \ldots = \tfrac{\beta_N}{\al^{(x)}_Na^2-\al^{(y)}_Nb^2}
\eeq
via the choice $\theta =- \tan^{-1}\left(\tfrac{\beta_j ab}{\al^{(x)}_ja^2-\al^{(y)}_jb^2}\right)$.

Note that in the case $N=2$ the condition \eqref{E:remove_cond} reduces to
\beq\label{E:remove_cond2}
a^2\left(\tfrac{\al^{(x)}_1}{\beta_1}-\tfrac{\al^{(x)}_2}{\beta_2}\right)=b^2\left(\tfrac{\al^{(y)}_1}{\beta_1}-\tfrac{\al^{(y)}_2}{\beta_2}\right),
\eeq
which is solvable always unless $\tfrac{\al^{(x)}_1}{\beta_1}-\tfrac{\al^{(x)}_2}{\beta_2}=0$ and $\tfrac{\al^{(y)}_1}{\beta_1}-\tfrac{\al^{(y)}_2}{\beta_2}\neq 0$ or vice versa and unless $\text{sign} \left(\tfrac{\al^{(x)}_1}{\beta_1}-\tfrac{\al^{(x)}_2}{\beta_2}\right) = -\text{sign}\left(\tfrac{\al^{(y)}_1}{\beta_1}-\tfrac{\al^{(y)}_2}{\beta_2}\right)$.

As the mixed derivative in \eqref{E:CME} cannot be removed in all cases, it is important to study PML for this system with $\beta_j\neq 0$.

\section{PML Derivation}\label{S:derive_PML}
Since the derivation of perfectly matched layers and their analysis are performed only for the linear part of the CNLS system \eqref{E:CME}, the analysis will be using merely the linear part of one scalar equation due
to the fact that the system is diagonal in its linear part. The linear problem at hand, thus, reads
\beq\label{E:lin_scalar}
\ri\pa_t u + (\al^{(x)}\pa_x^2+\al^{(y)}\pa_y^2+\beta\pa_x\pa_y)u=0, \qquad (x,y)\in \R^2, \ t\geq 0
\eeq
with $\al^{(x)}\al^{(y)}> \beta^2, \ \text{sign}(\al^{(x)})=\text{sign}(\al^{(y)})$. Note that the mixed derivative is not removed in \eqref{E:lin_scalar} because the constructed PML will be used in the coupled system \eqref{E:CME}, where the removal is not always possible as discussed in the previous section.

\begin{figure}[h!]
  \begin{center}
    \includegraphics[scale=0.52]{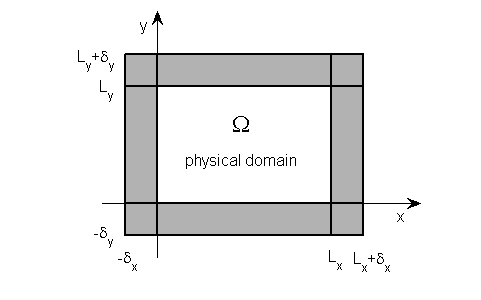}
  \end{center}
  \caption{Physical domain $\Omega = [0,L_x]\times [0,L_y]$ surrounded by layers.}
  \label{F:layer_sketch}
\end{figure}

Using the same approach as in \cite{DH07}, the PML is first derived in the directions of coordinate axes $x$ and $y$ and then combining the $x$ and $y-$layers, the corner layers are then proposed so that the resulting layer equations are applicable for all layers around the rectangular domain $\Omega$ as sketched in Fig. \ref{F:layer_sketch}. 

\begin{figure}[h!]
  \begin{center}
    \includegraphics[scale=0.45]{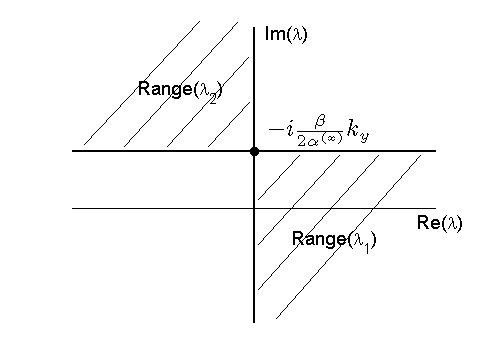}
  \end{center}
  \caption{Ranges of $\lambda_{1,2}$ for the modal solutions \eqref{E:modal_sol}.}
  \label{F:lambda_range}
\end{figure}

Without loss of generality the PML derivation is performed for the $x-$layer.
For that end let us suppose the domain is unbounded in $y$, so that $\Omega = [0,L_x]\times \R$. The Laplace transform in $t$ with $\text{Re}(s)\geq 0$ and Fourier transform in $y$ (with the dual variable $k_y$) of \eqref{E:lin_scalar} yield
\beq\label{E:transformed}
(\ri s + \al^{(x)}\pa_x^2-\al^{(y)}k_y^2+\ri\beta k_y\pa_x ) \hat{u} = 0.
\eeq
The initial data in the Laplace transform vanish because of the assumption $u(t=0)\equiv 0$ within the layers. The modal solutions of \eqref{E:transformed} are 
\beq\label{E:modal_sol}
\hat{u}(x;k_y,s) = e^{\lambda x}, \qquad \lambda = \lambda_{1,2}=\frac{1}{2\alpha^{(x)}}\left(-\ri\beta k_y\pm \sqrt{-\beta^2k_y^2-4\al^{(x)}(\ri s-\al^{(y)}k_y^2)}\right).
\eeq
The ranges of $\lambda_{1,2}$ for $\text{Re}(s)\geq 0$ are plotted in Fig. \ref{F:lambda_range}.

The analysis using Laplace transform is helpful as it allows for immediate identification of modes with positive and negative group velocity $v_g$. Indeed, for  $\text{Re}(s)\geq 0$ propagating modes with $v_g>0$ are contained within the modal set contained wholly in $\text{Re}(\lambda)\leq 0$, i.e., for the problem at hand in the $\lambda_2-$set; and those with $v_g<0$ are contained within the modal set contained wholly in $\text{Re}(\lambda)\geq 0$, i.e., in the $\lambda_1-$set. This can be seen by performing a Fourier transform in $x$ and Taylor expanding the dispersion $\omega(k_x)$ about $k_0\in \R$ for modes $e^{\ri(k_x x -\omega t)}$, where we set $\ri k_x = \lambda$ and $-\ri\omega =s$ \cite{Hag_priv_09}.

For \eqref{E:lin_scalar} the above relation between $\text{sign(Re}(\lambda))$ and $\text{sign}(v_g)$ can be easily checked by studying the dispersion relation explicitly. For the propagating modes $e^{\ri(k_x x -\omega t)}$ with $k_x,\omega\in \R$ the relation reads $\omega = \al^{(x)}k_x^2+\al^{(y)}k_y^2+\beta k_xk_y$ so that $v_g(k_x) = 2\al^{(x)}k_x +\beta k_y$. Clearly, $v_g > 0$ if and only if $k_x> {-\beta k_y \over 2\al^{(x)}}$. Let us stress that in the interval $k_x\in (0,{-\beta k_y \over 2\al^{(x)}})$ (or $k_x\in ({-\beta k_y \over 2\al^{(x)}},0)$ when ${\beta k_y \over 2\al^{(x)}}>0$) \textit{the group velocity and the phase velocity} $v_p(k_x) = \omega(k_x)/k_x$ \textit{have opposite signs}! This interval corresponds to the segment of the imaginary axis in Fig. \ref{F:lambda_range} between the origin and the point ${-\ri \beta k_y \over 2\al^{(x)}}$ since $\lambda = \ri k_x$. For many hyperbolic systems, like linearized Euler equations or equations for elastic waves, equality of sign$(v_p)$ and sign$(v_g)$ is shown to be a necessary condition for stability of the PML equations in time \cite{BFJ03,AK06}. Section \ref{S:stab} shows that in the Schr\"{o}dinger type problem at hand this is not the case and stability can be easily achieved by a choice of the PML parameters.

Let us now return to the transformed problem \eqref{E:transformed}. For absorption in the $x-$layers $x<0$ and $x>L_x$ the solution $\hat{u}$ needs to be modified (in the layers) to yield some $\uhatpml$ which in the layers satisfies
\beq\label{E:damp_prop}
\pa_x \uhatpml  = \tilde{\lambda} \uhatpml \quad \quad \text{with} \quad \left\{\begin{array}{rl}
\text{Re}(\tilde{\lambda}) <0 &\quad \text{for modes traveling right in} \ x\\
\text{Re}(\tilde{\lambda}) >0 &\quad \text{for modes traveling left in} \ x.
\end{array}\right.
\eeq
As the $\lambda$ in \eqref{E:modal_sol} satisfy the non-strict version of these inequalities, it would suffice to apply a simple rotation of $\lambda$ about $\lambda_0 = -\ri{ \beta k_y \over 2\al^{(x)}}$ by an angle $\rho\in (0,\pi/2)$, i.e., $\tilde{\lambda}=e^{\ri \rho}\left(\lambda+\ri{ \beta k_y \over 2\al^{(x)}}\right)-\ri{ \beta k_y \over 2\al^{(x)}}$. The resulting modes $\hat{u}^{\text{PML}} = e^{\tilde{\lambda}x}$ are, however, not perfectly matched with $\hat{u}$ at the interfaces $x=0$ and $x=L_x$. In order to achieve perfect matching one can, instead, set
\beq\label{E:u_pml}
\hat{u}^{\text{PML}} = e^{\left(\lambda +\ri{ \beta k_y \over 2\al^{(x)}}\right)\left(x+e^{\ri \rho}\int_{x_0}^x\sigma_x(\xi)d\xi\right)-\ri{ \beta k_y \over 2\al^{(x)}}x},
\eeq
where $x_0=0$ and $x_0=L_x$ for the layers $x<0$ and $x>L_x$ respectively. Choosing ${d^k\over dx^k}\sigma_x(x_0)=0$ for all $k\leq n-1$  guarantees $C^n$ matching of $\hat{u}$ and $\hat{u}^{\text{PML}}$ at $x=x_0$. With \eqref{E:u_pml} one obtains $\tilde{\lambda} = \lambda+\left(\lambda +\ri{ \beta k_y \over 2\al^{(x)}}\right)e^{\ri \rho}\sigma_x(x)$ in \eqref{E:damp_prop}. Under the condition $\sigma_x(x)>0$ for $x<0$ and $x>L_x$ this $\tilde{\lambda}$ can be seen using \eqref{E:modal_sol} or Fig. \ref{F:lambda_range} to satisfy the inequalities in \eqref{E:damp_prop}. Note that $\hat{u}^{\text{PML}}$ can be viewed as the solution over the whole domain $[-\delta_x,L_x+\delta_x]$ if for $x\in [0,L_x]$ one sets $\sigma_x(x)=0$.

To derive equations for $\upml$, let us express $\pa_x \hat{u}$ in terms of $\uhatpml$:
\[
\pa_x \hat{u} = e^{-\phi} {1\over 1+e^{\ri \rho}\sigma_x}\left(\pa_x-\ri k_y e^{\ri \rho}{ \beta \over 2\al^{(x)}} \sigma_x\right)\uhatpml,
\]
where $\phi = \left(\lambda +\ri{ \beta k_y \over 2\al^{(x)}}\right)e^{\ri \rho}\int_{x_0}^x\sigma_x(\xi)d\xi$. Thus, defining $\paxpml:={1\over 1+e^{\ri \rho}\sigma_x}\left(\pa_x-e^{\ri \rho}{ \beta \over 2\al^{(x)}}\sigma_x\pa_y \right)$ and eliminating the common factor $e^{-\phi}$, the $x-$layer equation reads
\beq\label{E:PML_equ_side}
\ri \pa_t\upml +\left(\al^{(x)}(\paxpml)^2 + \al^{(y)}\pa_y^2 +\beta \paxpml\pa_y\right)\upml=0.
\eeq

Treatment of the $y-$layers is completely analogous and defines the operator $\paypml:={1\over 1+e^{\ri \rho}\sigma_y}\left(\pa_y- e^{\ri \rho}{ \beta \over 2\al^{(y)}}\sigma_y\pa_x \right)$ with $\sigma_y=0$ for $y\in [0,L_y], \sigma_y(y)>0$ for $y<0$ and $y>L_y$ and with the perfect matching condition ${d^k\over dy^k}\sigma_y(y_0)=0$ for all $k\leq n-1$ at $y_0=0$ and $y_0=L_y$ with a chosen $n\in \N$. The most natural approach to the corner layers $[-\delta_x,0)\times[-\delta_y,0), \ [-\delta_x,0)\times(L_y,L_y+\delta_y], \ (L_x+\delta_x,L_x]\times (-\delta_y,0)$ and $(L_x+\delta_x,L_x]\times (L_y,L_y+\delta_y]$ is to combine the two layer equations into
\beq\label{E:PML_equ}
\ri \pa_t \upml +\left(\al^{(x)}(\paxpml)^2 + \al^{(y)}(\paypml)^2 +\beta \paxpml\paypml\right)\upml=0.
\eeq
Note that as 
\[
\left(\paxpml\paypml-\paypml\paxpml\right)\upml= \tfrac{e^{2\ri \rho}\beta }{ 2} \left[\tfrac{\sigma_x\sigma_y' \left({\beta \over 2\al^{(y)}}\pa_x+\pa_y\right) }{ \al^{(x)}(1+e^{\ri \rho}\sigma_y)}-\tfrac{\sigma_y\sigma_x' \left({\beta \over 2\al^{(x)}}\pa_y+\pa_x \right) }{ \al^{(y)}(1+e^{\ri \rho}\sigma_x)}\right],
\]
the operators $\paxpml$ and $\paypml$ do not commute unless $\beta=0$ or $\sigma_{x,y}$ are constant. Although in the performed numerical examples this does not seem to affect the $L^2$ error of the solution, in the numerical examples presented in Section \ref{S:num_tests} the operator $\paxpml\paypml$ was replaced by the commuting alternative ${1\over 2}\left(\paxpml\paypml+\paypml\paxpml\right)$.

In the nonlinear case $\Gamma \neq 0$ the above analysis does not, strictly speaking, apply. Nevertheless, if the solution remains small within the layers, the polynomial nonlinearity can be viewed as negligible and it makes sense to simply use
\beq\label{E:PML_equ_NL}
\ri \pa_t \upml_j +\left(\al_j^{(x)}(\paxpml)^2 + \al_j^{(y)}(\paypml)^2 +\beta_j \paxpml\paypml\right)\upml_j + \Gamma \CN_j(\upml_1, \ldots, \upml_N)=0, \quad  j \in \{1, \ldots, N\}
\eeq
as the PML system corresponding to \eqref{E:CME} in such a nonlinear scenario. This formulation is used in the nonlinear numerical simulations in Section \ref{S:num_tests_NL}. In fact, in one of the nonlinear numerical tests a large pulse enters the layer and the solution still qualitatively correct.

To the author's knowledge no truly perfectly matched layers for nonlinear systems exist in the literature. Appending the linear layer equations with the corresponding nonlinear terms is a common approach. For 
the nonlinear Schr\"{o}dinger equation (NLS) this is done, for instance, in \cite{FL05}. In \cite{Z07} the same type of PML is constructed non-rigorously by viewing the nonlinearity as a (solution dependent) potential. It is then argued that the success of this approach is due to the time-transverse invariant property of NLS. Radiation boundary conditions, on the other hand, have been successfully derived for some truly nonlinear systems including NLS \cite{S06,S06b}.

\section{Stability of the Layer Equations}\label{S:stab}

The analysis of Section \ref{S:derive_PML} does not guarantee that the layer equations are stable in time. Stability is determined below only for constant $\sigma_x$ and $\sigma_y$, in which case rewriting the layer equations in Fourier space results in a diagonal ODE system. The following analysis determines which parameter values (in particular $\sigma_{x,y}$) lead to boundedness of all Fourier modes in time, i.e., stability, and which lead to growth of at least one mode, i.e., instability. Once again, as the linear system is uncoupled, it is sufficient to study the scalar problem \eqref{E:PML_equ}. 

\subsection{Corner layer equations}\label{S:stab_corner}
\begin{figure}[h!]
  \begin{center}
    \includegraphics[scale=0.45]{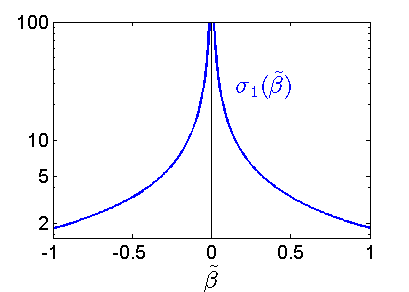}
  \end{center}
  \caption{The stability threshold function $\sigma_1(\tilde{\beta})$ in \eqref{E:D_roots}.}
  \label{F:sig_1}
\end{figure}

Let us firstly perform the change of variables $\tilde{x} = |\al^{(x)}|^{-1/2}x, \tilde{y} = |\al^{(y)}|^{-1/2}y$, which replaces $\al^{(x)}\rightsquigarrow 1, \al^{(y)}\rightsquigarrow 1$ and $\beta\rightsquigarrow |\al^{(x)}|^{-1/2}|\al^{(y)}|^{-1/2}\beta =:\tilde{\beta}$. Due to \eqref{E:coef_cond} one gets $|\tilde{\beta}|<1$. Dropping the tildes over the spatial variables, one obtains for the corner layer equations
\beq\label{E:PML_equ_scaled}
\ri \pa_t \upml +\left((\paxpml)^2 + (\paypml)^2 +\tilde{\beta} \paxpml\paypml\right)\upml=0.
\eeq
Applying the Fourier transform in $x$ and $y$, with dual variables $k_x$ and $k_y$, to \eqref{E:PML_equ_scaled} yields
\[
\begin{split}
\pa_t \uhatpml &= -\ri \left[k_x^2\left({1\over \mu_x^2}+{\tilde{\beta}^2e^{2\ri \rho}\sigma_y^2 \over 4\mu_y^2} - {\tilde{\beta}^2e^{\ri \rho}\sigma_y \over 2\mu_x\mu_y} \right)+k_y^2 \left({1\over \mu_y^2}+{\tilde{\beta}^2e^{2\ri \rho}\sigma_x^2 \over 4\mu_x^2} - {\tilde{\beta}^2e^{\ri \rho}\sigma_x \over 2\mu_x\mu_y} \right)\right.\\
&\quad \quad \quad\left. -\tilde{\beta}k_xk_y \left({e^{\ri \rho}\sigma_x \over \mu_x^2}+{e^{\ri \rho}\sigma_y \over \mu_y^2} - {4+e^{2\ri \rho}\tilde{\beta}^2\sigma_x\sigma_y\over 4\mu_x\mu_y} \right)\right]\uhatpml\\
&=:\nu(k_x,k_y)\uhatpml
\end{split}
\]
with $\mu_x = 1+e^{\ri \rho}\sigma_x$ and $\mu_y = 1+e^{\ri \rho}\sigma_y$. The stability requirement is $\text{Re}(\nu(k_x,k_y))\leq 0 \ \forall (k_x,k_y)\in \R^2$. For the sake of simplicity let us set $\sigma_x=\sigma_y=:\sigma$ and $\rho=\pi/4$. Under these simplifications 
\[
\begin{split}
\text{Re}(\nu) = {\sigma \over 4(\sigma^2+\sqrt{2}\sigma+1)^2}&\left[(k_x^2+k_y^2)\left(2\sqrt{2}\tilde{\beta}^2\sigma^2+\sigma(\tilde{\beta}^2-4)-\sqrt{2}(\tilde{\beta}^2+4)\right)\right.\\
&\ \left.+k_xk_y\tilde{\beta}\left(\sqrt{2}\sigma^2(\tilde{\beta}^2+4)+\sigma(\tilde{\beta}^2-4)-8\sqrt{2}\right)\right].
\end{split}
\]
Clearly, for $\tilde{\beta}\neq 0$ taking $\sigma>0$ small enough yields $\text{Re}(\nu)<0$ while large $\sigma>0$ result in $\text{Re}(\nu)>0$. Next, $\text{Re}(\nu)=0$ if and only if 
$$
k_x = k_y \tfrac{4-\tilde{\beta}^2}{2\tilde{\beta}^2(2\sigma^2-1)+ \sqrt{2}\sigma (\tilde{\beta}^2-4)-8} \left[\tilde{\beta}\left(\tfrac{\sigma}{\sqrt{2}} -\sigma^2\tfrac{\tilde{\beta}^2+4}{4-\tilde{\beta}^2}+\tfrac{8}{4-\tilde{\beta}^2}\right) \pm \left(\tilde{\beta}^2\sigma^4 +\sqrt{2}\sigma(\tilde{\beta}^2\sigma^2-4)+\left(\tfrac{\tilde{\beta}^2}{2}-2\right)\sigma^2-4\right)^{1/2}\right].
$$
Clearly, nonexistence of real solutions $(k_x,k_y)$ is equivalent to $D:=\tilde{\beta}^2\sigma^4 +\sqrt{2}\sigma(\tilde{\beta}^2\sigma^2-4)+\left(\tfrac{\tilde{\beta}^2}{2}-2\right)\sigma^2-4<0$. The roots of $D$ are 
\beq\label{E:D_roots}
\sigma_{1,2}(\tilde{\beta}) = \tfrac{\sqrt{2}}{4\tilde{\beta}}\left(2-\tilde{\beta}\pm \sqrt{\tilde{\beta}^2+12\tilde{\beta}+4}\right), \qquad \sigma_{3,4}(\tilde{\beta}) = -\tfrac{\sqrt{2}}{4\tilde{\beta}}\left(2+\tilde{\beta}\pm \sqrt{\tilde{\beta}^2-12\tilde{\beta}+4}\right).
\eeq
A straightforward analysis of $\sigma_{1,2,3,4}$ reveals that $D<0$ holds for $0< \tilde{\beta}<1$ when $0<\sigma<\sigma_1(\tilde{\beta})$ and for $-1< \tilde{\beta}<0$ when $0<\sigma<\sigma_3(\tilde{\beta})$. Because $\sigma_3(\tilde{\beta})=\sigma_1(-\tilde{\beta})$, this reduces to
\beq\label{E:stab_cond}
0<\sigma<\sigma_1(\tilde{\beta}) \quad \text{for} \quad |\tilde{\beta}|<1.
\eeq
Figure \ref{F:sig_1} plots the function $\sigma_1(\tilde{\beta})$. 

As a conclusion, under the condition \eqref{E:stab_cond} the layer equations \eqref{E:PML_equ} with constant $\sigma_x\equiv\sigma_y\equiv\sigma, \rho = \pi/4, \al^{(x)}=\al^{(y)}=1$ and $|\tilde{\beta}|<1$ are thus \textit{stable}. Note that $\sigma_1(\tilde{\beta})\rightarrow \infty$ as $\tilde{\beta}\rightarrow 0$ so that the layer equations are unconditionally stable for $\tilde{\beta}=0$, i.e., for the classical 2D Schr\"{o}dinger  equation. 

In order to determine the stability condition for the linear \textit{system} \eqref{E:PML_equ_NL} with $\Gamma=0$ each equation can be first scaled so that the coefficients of the non-mixed derivatives become $1$ and the mixed derivative coefficients become $\tilde{\beta}_j:=|\al_j^{(x)}|^{-1/2}|\al_j^{(y)}|^{-1/2}\beta_j$. The stability condition for $(\sigma_x)_j\equiv(\sigma_y)_j\equiv\sigma=\text{const.}$ is then
\beq\label{E:stab_cond_syst}
0<\sigma<\sigma_1\left(\max_{j\in \{1,\ldots,N\}}|\tilde{\beta}_j|\right).
\eeq
Because in practice $\sigma_{x,y}$ are usually taken non-constant, the above condition on $\sigma$ translates to a condition on $\max(\sigma_x)$ and $\max(\sigma_y)$. This can be justified by approximating $\sigma_{x,y}$ by piecewise constant functions and applying the stability condition \eqref{E:stab_cond_syst} on each piece.

Clearly, the presence of waves with group and phase velocities of opposite sign (see Section \ref{S:derive_PML}) does not lead to unconditional instability of the layer equations \eqref{E:PML_equ}, which is in contrast with the studies of several hyperbolic systems \cite{BFJ03,AK06}.

\subsection{Side layer equations}\label{S:stab_side}

Taking the Fourier transform of the side layer equation \eqref{E:PML_equ_side} under the assumption of $\sigma_x \equiv \sigma =$const. and defining $\nu$ analogously to Section \ref{S:stab_corner}, gives
\[\text{Re}(\nu) = -\tfrac{\sigma(\sqrt{2}+\sigma)}{(\sigma^2+\sqrt{2}\sigma+1)^2}(2k_x+\tilde{\beta} k_y)^2,\]
so that $\text{Re}(\nu)\leq 0 \quad \forall (k_x,k_y)\in \R^2$. 

Side layers are, therefore, unconditionally stable even for $\tilde{\beta} \neq 0$ and it is possible to use an absorption function $\sigma_x$, whose maximum in the corners satisfies \eqref{E:stab_cond_syst} and takes a larger value in the side layers $y\in [0,L_y]$ (and analogously for $\sigma_y$) so that the absorption in the side layers is strengthened. This is, however, not done in the simulations in Section \ref{S:num_tests} and $\sigma_x$ and $\sigma_y$ are kept $y$ and $x$ independent respectively. 

\section{Numerical Tests}\label{S:num_tests}

This section presents results of several numerical simulations of the system \eqref{E:PML_equ_NL} in both the linear ($\Gamma=0$) and the nonlinear ($\Gamma\neq 0$) case. 
The primary objective is to demonstrate convergence of the solution error with respect to the layer width. In the linear regime layers of infinite width ($\delta_x=\delta_y=\infty$) do not generate any error and the restrictions of the solution of \eqref{E:CME} and the solution of \eqref{E:PML_equ_NL} onto the physical domain $\Omega = [0,L_x]\times[0,L_y]$ are identical. Finite layers produce reflections from their far end but the resulting error in $\Omega$ decays exponentially with the layer width \cite{H99,deHoopTD,DJ06} due to the exponential decay of the solution within the layers. This exponential error convergence is numerically demonstrated to hold also here. Even in the nonlinear tests, where exponential convergence cannot be proved for the proposed layer equations, the resulting convergence is apparently exponential although the relative error in the example with large data in the layers becomes large.

In all the presented  numerical examples the rectangular domain $[-\delta_x,L_x+\delta_x]\times [-\delta_y,L_y+\delta_y]$ with $\delta_x,\delta_y,L_x,L_y>0$ was used with $\Omega=[0,L_x]\times[0,L_y]$ being the physical domain and the rest being the PML layers. The spatial discretization of the PDEs was done via the centered 4th order finite difference formulas 
\beq\label{E:FD_formulas}
\begin{split}
\pa_x^2u(x_i,y_j) &\approx (-u_{i-2,j}+16u_{i-1,j}-30u_{i,j}+16u_{i+1,j}-u_{i+2,j})/(12dx^2), \\
\pa_xu(x_i,y_j) &\approx (u_{i-2,j}-8u_{i-1,j}+8u_{i+1,j}-u_{i+2,j})/(12dx), 
\end{split}
\eeq
where $dx = x_{i+1}-x_i$ and $u_{i,j}=u(x_i,y_j)$; and analogously for the $y-$derivatives. The zero Dirichlet boundary condition was imposed at the outer layer boundary. The time-evolution was approximated via 4th order additive Runge-Kutta scheme of the ESDIRK type \cite{KC03}, in which the linear (stiff) terms are treated implicitly and the nonlinear terms explicitly.

Regarding the PML parameters, $\rho$ was taken $\rho=\pi/4$ and the absorption functions $\sigma_x$ and $\sigma_y$ were chosen of the form
\begin{equation}\label{E:sigma_fn}
\sigma_x(x)=\left\{ 
\begin{array}{ll}
\tfrac{h_x}{4}[1+\tanh(a_x(\delta_x)(x-L_x-\tfrac{\delta_x}{2}))][1+\tanh(6a_x(\delta_x)(x-L_x-\tfrac{\delta_x}{8}))] &\text{for} 
\quad x\in (L_x, L_x+\delta_x]\\
\tfrac{h_x}{4}[1-\tanh(a_x(\delta_x)(x+\tfrac{\delta_x}{2}))] [1-\tanh(6a_x(\delta_x)(x+\tfrac{\delta_x}{8}))] &\text{for} \quad
x\in [-\delta_x ,0)
\end{array}\right.
\end{equation}
with $a_x(\delta_x)=12/\delta_x$ and analogously for $\sigma_y(y)$. Clearly, $\max \sigma_x =  \tfrac{h_x}{4}[1+\tanh(a_x(\delta_x)\tfrac{\delta_x}{2})]\cdot$$[1+\tanh(a_x(\delta_x)\tfrac{21\delta_x}{4}))]$ is well approximated by $h_x$ even for moderate values of $\delta_x$. 
\begin{figure}[h!]
  \begin{center}
    \includegraphics[scale=0.45]{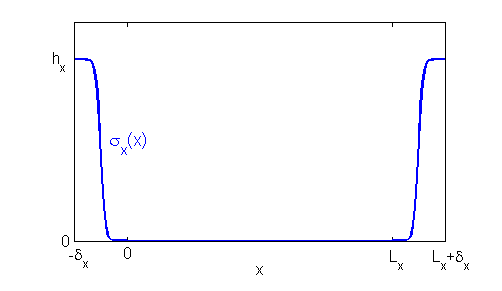}
  \end{center}
  \caption{The auxiliary PML function $\sigma_x(x)$ in \eqref{E:sigma_fn}.}
  \label{F:sig_plot}
\end{figure} 
The plot of $\sigma_x$ for $L_x/\delta_x=5$ is in Fig. \ref{F:sig_plot}. Note that the function $\sigma_x$ can certainly be chosen differently than \eqref{E:sigma_fn} and no claim on optimality is made here. The function in the second pair of square brackets on each line of \eqref{E:sigma_fn} is used merely to make $\sigma_x$ converge to $0$ at $x=0$ and $x=L_x$ in a smoother manner. One could, of course, drop this function and simply make $a_x(\delta_x)$ larger but that would result in a large slope of $\sigma_x$ within the layer, which leads to reflections in the numerical solution. An optimization study on the PML parameters, primarily $\sigma_x$ and $\sigma_y$, can be performed \cite{CM98} to increase efficiency of the layers.

\subsection{Simulations of the Linear Case $\Gamma = 0$}\label{S:num_tests_lin}

Clearly, in the linear case $\Gamma=0$ the system \eqref{E:PML_equ_NL} decouples and the change of variables \eqref{E:transf} (possibly distinct for each $j$) can be applied to remove the cross-derivatives. Nevertheless, because in the to-be-studied nonlinear case this removal is not always possible (see Section \ref{S:CNLS_relevance}), one of the two linear numerical tests provided below is with the mixed derivatives present. Initial data localized at the center of $\Omega$ were used in the tests.
The error was computed at $t=1$ using the reference solution $\uref$ determined by numerically solving the initial value problem \eqref{E:PML_equ_NL} with $N=1, \Gamma=0$ via Fourier transform  on a much larger domain (than $\Omega$) on which the dispersed solution at $t=1$ is well localized. 
 
\subsubsection{Linear Scalar Case with $\beta=0$}\label{S:lin_beta_0}

In the linear scalar case $\Gamma=0, N=1$ with $\beta=0$ the problem reduces to the 2D linear Schr\"{o}dinger equation and the layer equations \eqref{E:PML_equ} are those used extensively in the literature, see e.g. \cite{H03_b,FL05,Z07}. The numerical test for this case is presented here for completeness and comparison with the case $\beta \neq 0$ as well as with other publications.

The remaining coefficients are chosen $\al^{(x)}=3/4$ and $\al^{(y)}=5/4$ and the numerical parameters are $L_x=L_y=6, dx=dy = L_x/350 \approx 0.017$ and $dt = 0.01$. The PML parameters are $\rho=\pi/4, h_x = 30$ and the computations were performed for 6 different layer widths $\delta_x=\delta_y \in \{0.08 L_x,0.12 L_x, 0.16 L_x, 0.2 L_x, 0.25 L_x, 0.3 L_x\}$ with the initial data $u(x,y,0) = e^{-(x-L_x/2)^2-(y-L_y/2)^2}$.
\begin{figure}[h!]
  \begin{center}
    \includegraphics[scale=0.45]{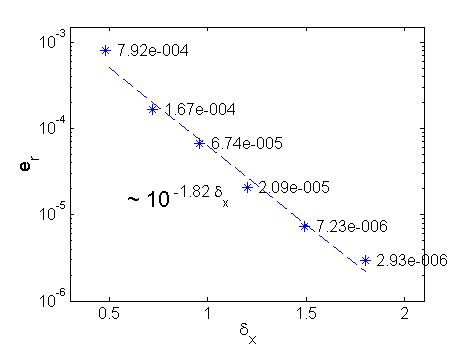}
  \end{center}
  \caption{Error convergence with respect to the layer width for the test in Section \ref{S:lin_beta_0}. stars: relative $L^2$ error $e_r = \|\upml - \uref\|_{L^2(\Omega)}/\|\uref\|_{L^2(\Omega)}$ at $t=1$, dashed line: $c \ 10^{-1.82 \ \delta_x}$.}
  \label{F:err_conv_lin_1}
\end{figure} 
Fig. \ref{F:err_conv_lin_1} shows convergence of the $L^2$ norm of the error over the physical domain $\Omega$ at $t=1$ in dependence on $\delta_x=\delta_y$ featuring an exponential convergence $e^{-p\delta_x}$ with $p\approx 1.82$.

\subsubsection{Linear Scalar Case with $\beta=0.5$}\label{S:lin_beta_nzero}

The same parameters and initial data as in Section \ref{S:lin_beta_0} are chosen here except for the following: $\al^{(x)}=\al^{(y)}=1, \beta =0.5$ and the maximum value of $\sigma_{x,y}$, which is set to $h_x=h_y=3.3$, i.e., close to the stability threshold $\sigma_1(0.5)\approx 3.325$ in \eqref{E:stab_cond}.

\begin{figure}[h!]
  \begin{center}
    \includegraphics[scale=0.45]{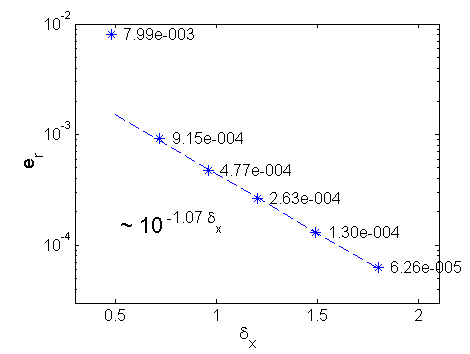}
  \end{center}
  \caption{Error convergence with respect to the layer width for the test in Section \ref{S:lin_beta_nzero}. stars: relative $L^2$ error $e_r = \|\upml - \uref\|_{L^2(\Omega)}/\|\uref\|_{L^2(\Omega)}$  at $t=1$, dashed line: $c \ 10^{-1.07 \ \delta_x}$.}
  \label{F:err_conv_lin_2}
\end{figure}

\begin{figure}[h!]
  \begin{center}
    \includegraphics[scale=0.45]{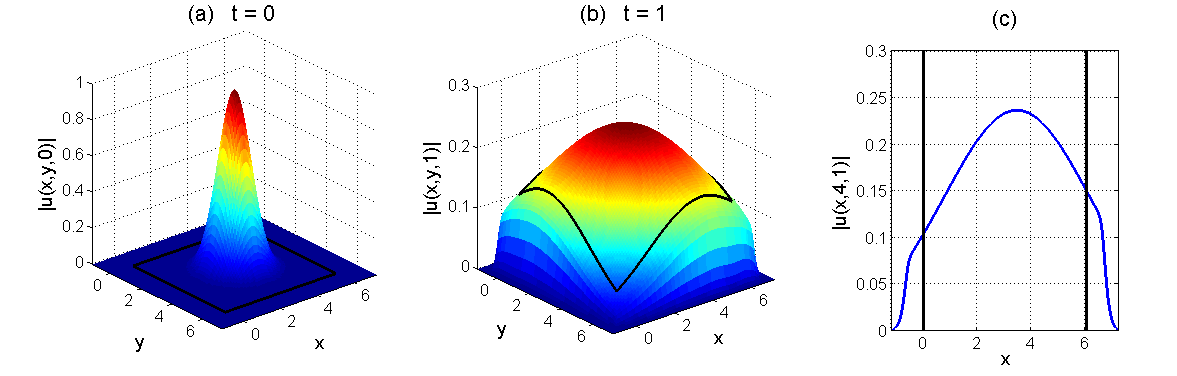}
  \end{center}
  \caption{The solution modulus for the test in Section \ref{S:lin_beta_nzero} with $\del_x=\del_y=0.2 L_x$. Black lines denote the interface between the physical domain and PML layers. (a) initial data; (b) solution modulus at $t=1$; (c) solution modulus at $t=1,y=4$.}
  \label{F:u_prof_lin_2}
\end{figure}

\begin{figure}[h!]
  \begin{center}
    \includegraphics[scale=0.45]{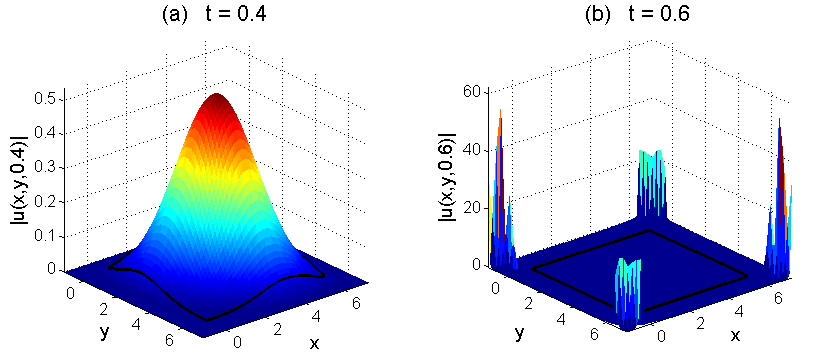}
  \end{center}
  \caption{The solution modulus for the unstable test in Section \ref{S:lin_beta_nzero} with $h_x=h_y = 20$. (a) solution modulus at $t=0.4$; (b) solution modulus at $t=0.6$.}
  \label{F:u_prof_unstable_PML}
\end{figure}

As one can see in Fig. \ref{F:err_conv_lin_2}, the error convergence is, once again, exponential like $e^{-p\delta_x}$ with $p\approx 1.07$. Compared to the $\beta=0$ case in Section \ref{S:lin_beta_0} the convergence is slower and the error values are slightly larger which is to be expected due to the weaker applied absorption.
Fig. \ref{F:u_prof_lin_2} shows the initial profile and the modulus of the solution at $t=1$ for $\delta_x=\delta_y = 0.2 L_x = 1.2$. A simulation with the same coefficients and similar discretization and PML parameters to those in Fig. \eqref{F:u_prof_lin_2} but with $h_x=h_y=20$, which exceeds the stability threshold \eqref{E:stab_cond}, is shown in Fig. \ref{F:u_prof_unstable_PML}. It, indeed, results in an instability within the layers, clearly seen in the plot at $t=0.6$.

\subsection{Simulations of the Nonlinear Case $\Gamma \neq 0$}\label{S:num_tests_NL}

As advertised in Section \ref{S:CNLS_relevance}, the CNLS system \eqref{E:CME2} was used for numerical tests in the nonlinear case and $\eps_q$ was taken negative in order to prevent possible finite time blowup of the solution \cite{SS00}. Three tests are performed below. In the first two tests (Sections \ref{S:NL_beta_0} and \ref{S:NL_beta_nzero}) the following choice of coefficients and numerical parameters was made: $\Gamma=0.5, \eps_q=-0.2, L_x=L_y=14, dx=dy = 14/250 \approx 0.056, dt=0.01, \rho =\pi/4 $ and 
$\del_x=\del_y\in \{0.08 L_x, 0.12 L_x, 0.16 L_x, 0.2 L_x, 0.25 L_x, 0.3 L_x\}.$  Section \ref{S:NL_beta_0} presents a case where the choice of $\al^{(x)}_{1,2},\al^{(y)}_{1,2}$ and $\beta_{1,2}$ allows a change of variables that removes the cross-derivatives and section \ref{S:NL_beta_nzero} a case where this is impossible. In both cases the initial data are the sum of a stationary solitary wave and four perturbing Gaussians so that the dynamics result in a large amount of radiation shed toward the boundary with a solitary waves remaining at the domain center. The solution is evolved up to $t=5$. In detail, the initial data are
\beq\label{E:NL_IC}
\vec{u}(x,y,0) = \vec{\phi}_s(x-L_x/2,y-L_y/2)+0.8\sum_{k=1}^4 e^{-2\left((x-p_k)^2+(y-q_k)^2\right)}\bpm 1\\1 \epm 
\eeq
with $p_{1,2,3,4} = \tfrac{L_x}{2},\tfrac{L_x}{2},\tfrac{L_x}{4},\tfrac{3L_x}{4}$ and $q_{1,2,3,4} = \tfrac{L_y}{4},\tfrac{3L_y}{4},\tfrac{L_y}{2},\tfrac{L_y}{2}$ respectively, and where $\vec{\phi}_s(x,y)$ is the positive spatial profile of the stationary solitary wave (ground state) $e^{-\ri t}  \vec{\phi}_s(x,y)$ of \eqref{E:CME2}. $\vec{\phi}_s(x,y)$ was computed via Newton's iteration on the corresponding stationary system, i.e., on \eqref{E:CME2} with $\ri \pa_t$ replaced by $1$ and with zero Dirichlet boundary conditions at $\pa \Omega$. In detail \eqref{E:CME2} is first solved with $\beta_1=\beta_2=0$ for the radially symmetric Townes soliton with $u_1=u_2$ via the shooting method. Next, the solution is numerically continued via homotopy in $\beta_1$ and $\beta_2$ solving for $u_1$ and $u_2$ via Newton's iteration at each $\beta$-step and using the previous solution as an initial guess.

The third example (Section \ref{S:NL_pulse}) tests the designed PML for the scenario of a pulse entering the layers in the nonlinear regime; for the choice of parameters see the corresponding Section.

For all three tests exponential convergence of the error is observed despite the fact that the problem is nonlinear and there is no guarantee for such a convergence. In the first two tests the relative error is satisfactory while in the third case  which is truly nonlinear even in the layers, the relative error is rather large but the solution is still qualitatively correct. Note that for all three tests below the figures with the solution profiles show only the first component $u_1$ as $u_2$ behaves in qualitatively the same way.

\subsubsection{Nonlinear System with $\beta_j=0$}\label{S:NL_beta_0}
 The coefficients here are $\al^{(x)}_1 =\al^{(y)}_2= \tfrac{3}{4},\al^{(x)}_2 = \al^{(y)}_1=\tfrac{5}{4}$ and $\beta_1=\beta_2=0$, which can be viewed as obtained from the system with $\al^{(x)}_1 =\al^{(x)}_2 =\al^{(y)}_1 =\al^{(y)}_2 =1$ and $\beta_1=\beta_2=\tfrac{1}{2}$ via the transformation \eqref{E:transf} with $a=b=1$ and $\theta=\tfrac{\pi}{4}$. The magnitude of the absorption functions $\sigma_{x,y}$ is $h_x=h_y=8$.
 
\begin{figure}[h!]
  \begin{center}
    \includegraphics[scale=0.45]{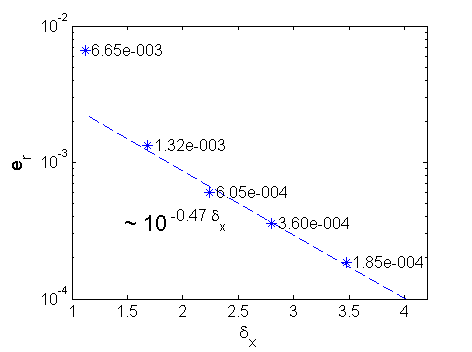}
  \end{center}
  \caption{Error convergence with respect to the layer width for the test in Section \ref{S:NL_beta_0}. stars: relative $L^2$ error $e_r = \|\upml - \uref\|_{L^2(\Omega)}/\|\uref\|_{L^2(\Omega)}$ at $t=5$, dashed line: $c \ 10^{-0.47 \ \delta_x}$.}
  \label{F:err_conv_NL_1}
\end{figure}

\begin{figure}[h!]
  \begin{center}
    \includegraphics[scale=0.45]{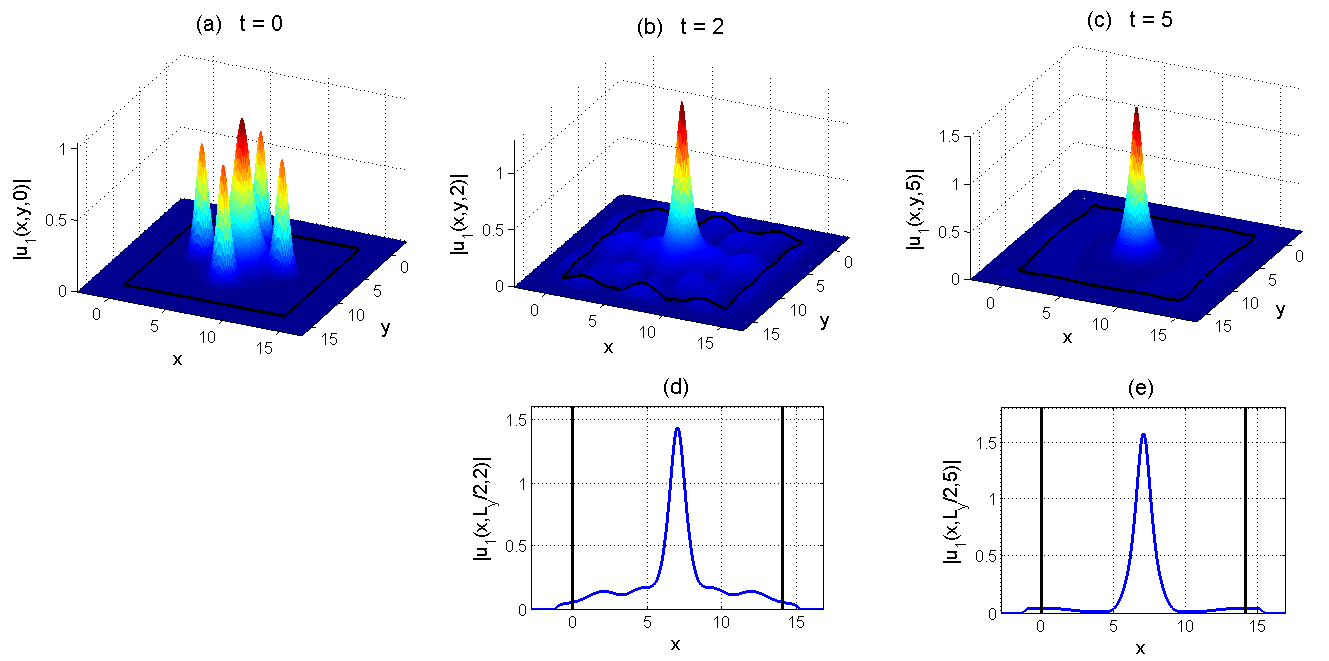}
  \end{center}
  \caption{The solution modulus for the test in Section \ref{S:NL_beta_0} with $\del_x=\del_y=0.2 L_x$. Black lines denote the interface between the physical domain and PML layers. (a) initial data; (b) solution modulus at $t=2$; (c) solution modulus at $t=5$; (d) and (e) solution modulus along $y=L_y/2$ at $t=2$ and $t=5$ respectively.}
  \label{F:u_prof_NL_1}
\end{figure}

 Fig. \ref{F:err_conv_NL_1} presents the error convergence at $t=5$ with respect to the layer width, where the solution computed with the widest layer ($\del_x=\del_y=0.3 L_x$) was used as the reference solution $\uref$. The convergence is exponential, like $e^{-p\delta_x}$ with $p\approx 0.47$.
Fig. \ref{F:u_prof_NL_1} shows for $\del_x=\del_y=0.2 L_x=2.8$ the initial data and the modulus of the first component $u_1$ at times $t=2$, when a large amount of radiation is traveling into the layers and at $t=5$, when most radiation has been absorbed.

\subsubsection{Nonlinear System with $\beta_j\neq0$}\label{S:NL_beta_nzero}

This example presents the case $\al^{(x)}_1 = 1,\al^{(x)}_2 = \tfrac{3}{4}, \al^{(y)}_1 = \al^{(y)}_2 = 1$ and $\beta_1=0.2, \beta_2=0.15$, for which equation \eqref{E:remove_cond2} cannot be solved for nonzero $a,b$ and thus the mixed derivatives cannot be removed. The initial condition \eqref{E:NL_IC} is, once again, used, where the solitary wave profile $\vec{\phi}_s$ differs from that used in Section \ref{S:NL_beta_0} due to the different PDE coefficients. The stability threshold for the selected coefficients as given by \eqref{E:stab_cond_syst} is $\sigma_1(\max(0.2,0.15\tfrac{2}{\sqrt{3}}))=\sigma_1(0.2)\approx 7.67$. In the simulation $h_x=h_y=7.6$ was used. 

\begin{figure}[h!]
  \begin{center}
    \includegraphics[scale=0.45]{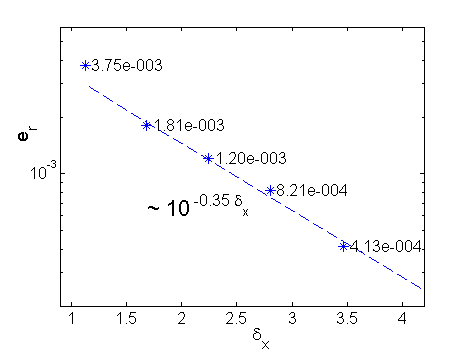}
  \end{center}
  \caption{Error convergence with respect to the layer width for the test in Section \ref{S:NL_beta_nzero}. stars: relative $L^2$ error $e_r = \|\upml - \uref\|_{L^2(\Omega)}/\|\uref\|_{L^2(\Omega)}$ at $t=5$, dashed line: $c \ 10^{-0.35 \ \delta_x}$.}
  \label{F:err_conv_NL_2}
\end{figure}

\begin{figure}[h!]
  \begin{center}
    \includegraphics[scale=0.45]{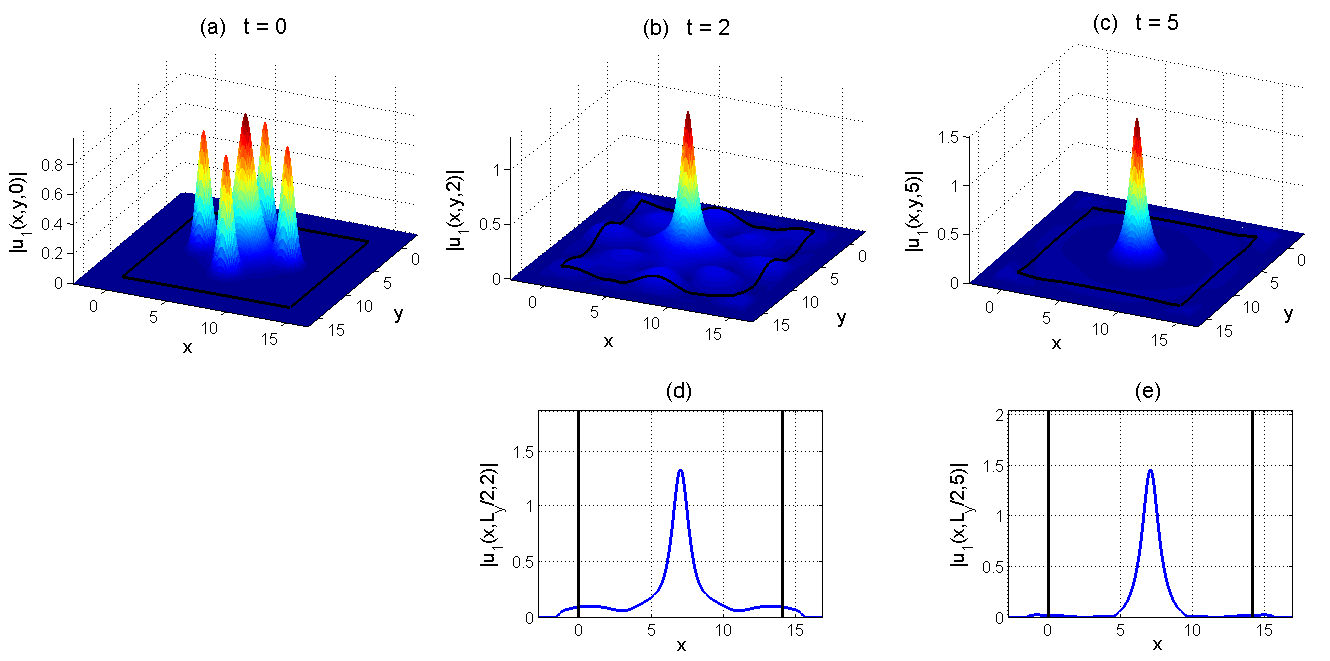}
  \end{center}
  \caption{The solution modulus for the test in Section \ref{S:NL_beta_nzero} with $\del_x=\del_y=0.2 L_x$. Black lines denote the interface between the physical domain and PML layers. (a) initial data; (b) solution modulus at $t=2$; (c) solution modulus at $t=5$; (d) and (e) solution modulus along $y=L_y/2$ at $t=2$ and $t=5$ respectively.}
  \label{F:u_prof_NL_2}
\end{figure}

Using, once again, the solution computed with the widest layer ($\del_x=\del_y=0.3 L_x$) as the reference solution, the error convergence is plotted in Fig. \ref{F:err_conv_NL_2}, where exponential convergence $e^{-p\del_x}, p \approx 0.35$ can be observed. The convergence rate is slightly smaller than that in Section \ref{S:NL_beta_0} due to the weaker applied absorption. Fig. \ref{F:u_prof_NL_2} then shows the solution modulus at selected instances of time.
Fig. \ref{F:u_prof_NL_2} shows for $\del_x=\del_y=0.2 L_x=2.8$ the initial data and the modulus of the first component $u_1$ at times $t=2$, when a large amount of radiation is traveling into the layers and $t=5$ when most radiation has been absorbed.

\begin{figure}[h!]
  \begin{center}
    \includegraphics[scale=0.45]{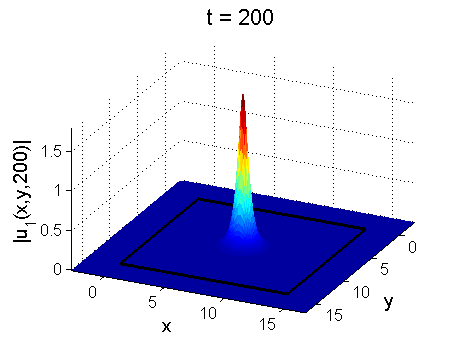}
  \end{center}
  \caption{The solution from Fig. \ref{F:u_prof_NL_2} at $t=200$.}
  \label{F:u_prof_NL_2_long_t}
\end{figure}

Finally, to check the stability result in long time dynamics, Fig. \ref{F:u_prof_NL_2_long_t} shows the solution profile at $t=200$ and demonstrates that no growth occurs.

\subsubsection{Nonlinear System with $\beta_j\neq0$ and a Pulse Propagating into the Layer}\label{S:NL_pulse}

\begin{figure}[h!]
  \begin{center}
    \includegraphics[scale=0.45]{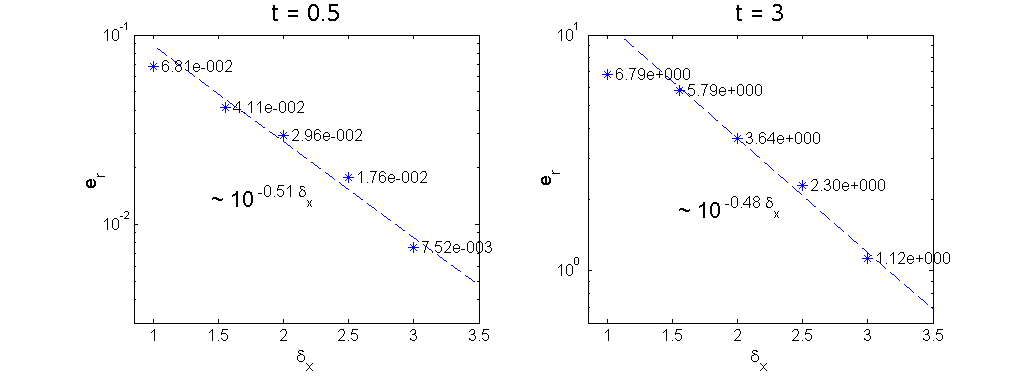}
  \end{center}
  \caption{Error convergence with respect to the layer width for the test in Section \ref{S:NL_pulse} on the left error at $t=0.5$ and on the right at $t=3$. stars: relative $L^2$ error $e_r = \|\upml - \uref\|_{L^2(\Omega)}/\|\uref\|_{L^2(\Omega)}$, dashed line: $c \ 10^{-0.51 \ \delta_x}$ and $c \ 10^{-0.48 \ \delta_x}$ on the left and right respectively.}
  \label{F:err_conv_pulse}
\end{figure}

\begin{figure}[h!]
  \begin{center}
    \includegraphics[scale=0.47]{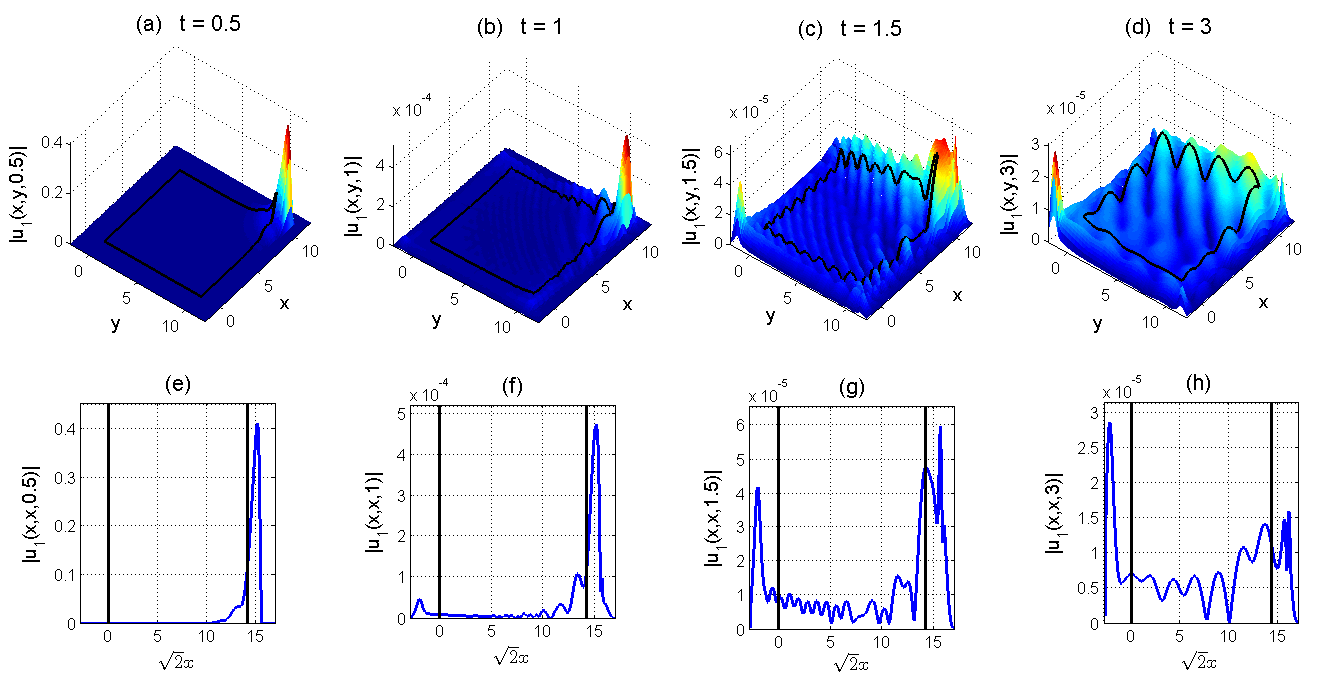}
  \end{center}
  \caption{The solution modulus for the test in Section \ref{S:NL_pulse} with $\delta_x=\delta_y= 0.2 L_x$. Black lines denote the interface between the physical domain and PML layers. (a) - (d) solution modulus at $t=0.5, 1, 1.5$ and $t=3$ resp.; (e) - (h) solution modulus along $y=x$ corresponding to (a)-(d) respectively.}
  \label{F:u_prof_pulse}
\end{figure}

In order to test the performance of the proposed layer equations in a truly nonlinear regime, Figs. \ref{F:err_conv_pulse} and \ref{F:u_prof_pulse} present the case of a pulse propagating into the layer (entering it at a corner of $\Omega$). The PDE coefficients are taken the same as in Sec. \ref{S:NL_beta_nzero} and the initial data used were the stationary solitary wave used in Sec. \ref{S:NL_beta_nzero}, centered at $(x,y)=(L_x/2,L_y/2)$ and multiplied by a plane wave in order to induce motion of the pulse along the $y=x$ line
\beq\label{E:NL_IC_kick}
\vec{u}(x,y,0) = \vec{\phi}_s(x-L_x/2,y-L_y/2) e^{6i\left((x-L_x/2)+(y-L_y/2)\right)}. 
\eeq
This initial condition does not correspond to an exact moving solitary wave solution of \eqref{E:CME2} (such solutions have not been found for this system with $\beta_1\neq \beta_2$) but the resulting solution propagates in a close to solitary manner.

The physical domain $\Omega=[0,L_x]\times[0,L_y]$ was set by $L_x=L_y=10$, the spatial discretization by $dx=dy = L_x/180\approx 0.056$ and the simulation was carried out for six different layer widths 
$\delta_x=\delta_y \in \{0.1 L_x,0.15 L_x, 0.2 L_x,$ $0.25 L_x, 0.3 L_x, 0.35 L_x\}$. The rest of the PML parameters was as in Section \ref{S:NL_beta_nzero}. Using the solution with $\delta_x=0.35 L_x$ as the reference solution, the convergence of the relative $L^2$ error over $\Omega$ is plotted in Fig. \ref{F:err_conv_pulse} at both $t=0.5$ when the pulse has just entered the layer and at $t=3$ when the pulse has propagated far from the domain $\Omega$ and the exact solution in $\Omega$ consists only of an exponentially small tail plus small radiation due to the fact that the pulse is not an exact traveling wave. 
The relative error is much larger than in the previous examples, which were linear or effectively linear in the layers, but on the selected $\delta_x$ range the convergence seems to be again exponential at both $t=0.5$ and $t=3$. More importantly, as one can see in Fig. \ref{F:u_prof_pulse}, the qualitative behavior of the solution is correctly captured and the pulse, which has amplitude about $0.99$ before entering the layers, leaves $\Omega$ with only very small reflected waves (amplitude $\sim 10^{-5}$) remaining. 
Fig. \ref{F:u_prof_pulse} shows the solution modulus at several instances of time over the whole spatial domain as well as along the line $y=x$, along which the pulse propagates. The first shown instance is at $t=0.5$ because for $t<0.5$ the pulse is simply traveling from its initial location at $(L_x/2,L_y/2)$ toward the corner layer. 

The results of this test suggest that the layer equations can be applied even in many truly nonlinear cases, where the solution amplitude within the layers may become large, mainly if only qualitative behavior of the pulses is required. No guarantee can, however, be given that its performance will be satisfactory in all such cases. Examples of relevant nonlinear problems are interaction of several pulses or interaction of pulses with localized defects, where one or more pulses leaves $\Omega$ within the simulation time.

\section{Discussion}\label{S:discussion}
The presented analysis of PML for the 2D Schr\"{o}dinger equation with cross derivatives shows that the presence of the cross derivatives leads to the existence of linear (Fourier) modes with opposite group and phase velocities. Unlike in some hyperbolic systems \cite{BFJ03,AK06} the resulting layer equations are only conditionally unstable and a choice of the damping functions $\sigma_x,\sigma_y$ below a calculatable threshold leads to stability of the linear PML equations. Note that the damping of the PML is ensured based on an analysis of group velocity rather than phase velocity of the linear modes. 

In the nonlinear case the linear PML equations are simply appended with the (polynomial) nonlinear terms. The layer performance is then affected only slightly if the solution remains small in the layers as seen in the provided numerical tests. An analysis of the layer performance in the nonlinear case would be valuable. Of tremendous interest would a perfectly matched layer for truly nonlinear waves, i.e. without the smallness assumption. Such analysis does not appear in the literature. In the area of radiation boundary conditions, on the other hand, limited results for nonlinear equations exist, see \cite{S06,S06b}.

The paper studies PML in the 2D case. Nevertheless, in 3D the derivation is completely analogous. The linear scalar equation corresponding to \eqref{E:lin_scalar} is in 3D
\beq\label{E:lin_scalar_3D}
\ri\pa_t u + (\al^{(x)}\pa_x^2+\al^{(y)}\pa_y^2+\al^{(z)}\pa_z^2+\beta_1\pa_x\pa_y+\beta_2\pa_x\pa_z+\beta_3\pa_y\pa_z)u=0, \qquad (x,y,z)\in \R^3, \ t\geq 0.
\eeq
The modal solutions analogous to \eqref{E:modal_sol} are $\hat{u}(x;k_y,k_z,s) = e^{\lambda x}$ with
\[
\qquad \lambda = \lambda_{1,2}=\frac{1}{2\alpha^{(x)}}\left(-\ri(\beta_1 k_y+\beta_2k_z)\pm \sqrt{-(\beta_1k_y+\beta_2k_z)^2-4\al^{(x)}(\ri s-\al^{(y)}k_y^2-\al^{(z)}k_z^2-\beta_3k_y k_z)}\right)
\]
and, thus, $\paxpml$ generalizes to
$\paxpml:={1\over 1+e^{\ri \rho}\sigma_x}\left(\pa_x-{e^{\ri \rho}\sigma_x \over 2\al^{(x)}}(\beta_1\pa_y+\beta_2\pa_z) \right)$. The operators $\paypml$ and $\pazpml$ are defined analogously and the layer equations are similarly to \eqref{E:PML_equ} and \eqref{E:PML_equ_NL} obtained by replacing $\pa_x, \pa_y$ and $\pa_z$ by 
$\paxpml$, $\paypml$ and $\pazpml$ respectively.

The algebra in the stability analysis becomes in 3D, however, much more complicated and will be left for future work. Note that PML for 3D Schr\"{o}dinger equations with mixed derivatives have been previously used in \cite{Cheng_etal07}. Perfect matching and stability were, however, not analyzed there in the presence of mixed derivatives.

\bigskip

\noi
{\bf Acknowledgement.}  The work of T. Dohnal is supported by the Humboldt
Research Fellowship. The author wishes to thank Thomas Hagstrom, Southern Methodist University, for stimulating discussions and many helpful remarks.

\bibliographystyle{is-unsrt}
\bibliography{bibliography}

\begin{thebibliography}{10}

\bibitem{Beren94}
J.-P. B\'erenger.
\newblock A perfectly matched layer for the absorption of electromagnetic
  waves.
\newblock {\em J. Comput. Phys.}, 114\penalty0 (2):\penalty0 185 -- 200, 1994.

\bibitem{TY98}
E.~Turkel and A.~Yefet.
\newblock Absorbing {PML} boundary layers for wave-like equations.
\newblock {\em Appl. Numer. Math.}, 27\penalty0 (4):\penalty0 533 -- 557, 1998.
\newblock Special Issue on Absorbing Boundary Conditions.

\bibitem{H98}
J.~S. Hesthaven.
\newblock On the analysis and construction of perfectly matched layers for the
  linearized {E}uler equations.
\newblock {\em J. Comput. Phys.}, 142\penalty0 (1):\penalty0 129 -- 147, 1998.

\bibitem{H03}
T.~Hagstrom.
\newblock A new construction of perfectly matched layers for hyperbolic systems
  with applications to the linearized {E}uler equations.
\newblock In G.~Cohen, E.~Heikkola, P.~Joly, and P.~Neittaanm\"{a}ki, editors,
  {\em {Mathematical and Numerical Aspects of Wave Propagation Phenomena}},
  pages 125--129. Springer-Verlag, Berlin, 2003.

\bibitem{AHK06}
D.~Appel\"{o}, T.~Hagstrom, and G.~Kreiss.
\newblock Perfectly matched layers for hyperbolic systems: {G}eneral
  formulation, well-posedness and stability.
\newblock {\em SIAM J. Appl. Math.}, 67:\penalty0 1--23, 2006.

\bibitem{C97}
F.~Collino.
\newblock Perfectly matched absorbing layers for the paraxial equations,.
\newblock {\em J. Comput. Phys.}, 131\penalty0 (1):\penalty0 164 -- 180, 1997.

\bibitem{H03_b}
T.~Hagstrom.
\newblock New results on absorbing layers and radiation boundary conditions.
\newblock In M.~Ainsworth, P.~Davies, D.~Duncan, P.~Martin, and B.~Rynne,
  editors, {\em {Topics in Computational Wave Propagation}}, pages 1--42.
  Springer-Verlag, 2003.

\bibitem{DU09}
T.~Dohnal and H.~Uecker.
\newblock Coupled mode equations and gap solitons for the 2d
  {G}ross-{P}itaevskii equation with a non-separable periodic potential.
\newblock {\em Physica D}, 238\penalty0 (9-10):\penalty0 860 -- 879, 2009.

\bibitem{DH07}
T.~Dohnal and T.~Hagstrom.
\newblock Perfectly matched layers in photonics computations: 1d and 2d
  nonlinear coupled mode equations.
\newblock {\em J. Comput. Phys.}, 223\penalty0 (2):\penalty0 690 -- 710, 2007.

\bibitem{BEVK05}
A.D. Boardman, P.~Egan, L.~Velasco, and N.~King.
\newblock Control of planar nonlinear guided waves and spatial solitons with a
  left-handed medium.
\newblock {\em J. Opt. A: Pure Appl. Opt.}, 7\penalty0 (2):\penalty0 S57--S67,
  2005.

\bibitem{DA05}
T.~Dohnal and A.B. Aceves.
\newblock {Optical soliton bullets in (2+1)D nonlinear Bragg resonant periodic
  geometries}.
\newblock In J.~Yang, editor, {\em {Nonlinear Wave Phenomena in Periodic
  Photonic Structures}}, volume 115 of {\em {Studies in Appl. Math.}}, pages
  209--232. 2005.

\bibitem{GW08}
R.H. Goodman and M.I. Weinstein.
\newblock Stability and instability of nonlinear defect states in the coupled
  mode equations—analytical and numerical study.
\newblock {\em Physica D}, 237:\penalty0 2731--2760, 2008.

\bibitem{BFJ03}
E.~B\'ecache, S.~Fauqueux, and P.~Joly.
\newblock Stability of perfectly matched layers, group velocities and
  anisotropic waves.
\newblock {\em J. Comput. Phys.}, 188\penalty0 (2):\penalty0 399 -- 433, 2003.

\bibitem{AK06}
D.~Appel\"{o} and G.~Kreiss.
\newblock A new absorbing layer for elastic waves.
\newblock {\em J. Comput. Phys.}, 215\penalty0 (2):\penalty0 642 -- 660, 2006.

\bibitem{Cheng_etal07}
C.~Cheng, J.-H. Lee, K.H. Lim, H.Z. Massoud, and Q.H. Liu.
\newblock 3d quantum transport solver based on the perfectly matched layer and
  spectral element methods for the simulation of semiconductor nanodevices.
\newblock {\em J. Comp. Phys.}, 227\penalty0 (1):\penalty0 455 -- 471, 2007.

\bibitem{SY07}
Z.~Shi and J.~Yang.
\newblock {Solitary waves bifurcated from Bloch-band edges in two-dimensional
  periodic media}.
\newblock {\em Phys. Rev. E}, 75\penalty0 (5):\penalty0 056602, 2007.

\bibitem{DPS09}
T.~Dohnal, D.~Pelinovsky, and G.~Schneider.
\newblock Coupled-mode equations and gap solitons in a two-dimensional
  nonlinear elliptic problem with a separable periodic potential.
\newblock {\em {J. Nonlin. Sci.}}, 19:\penalty0 95--131, 2009.

\bibitem{SS00}
C.~Sulem and P.~Sulem.
\newblock {\em {The nonlinear Schr\"{o}dinger equation: self-focusing and wave
  collapse}}.
\newblock Springer-Verlag, Berlin, 2000.
\newblock 93-103 pp.

\bibitem{Hag_priv_09}
T.~Hagstrom, 2009.
\newblock Private Discussions.

\bibitem{FL05}
C.~Farrell and U.~Leonhardt.
\newblock {The perfectly matched layer in numerical simulations of nonlinear
  and matter waves}.
\newblock {\em {J. Opt. B: Quantum Semiclass. Opt.}}, 7:\penalty0 1 -- 4, 2005.

\bibitem{Z07}
C.~Zheng.
\newblock A perfectly matched layer approach to the nonlinear {S}chr\"{o}dinger
  wave equations.
\newblock {\em J. Comput. Phys}, 227\penalty0 (1):\penalty0 537 -- 556, 2007.

\bibitem{S06}
J.~Szeftel.
\newblock Absorbing boundary conditions for one-dimensional nonlinear
  {S}chr\"odinger equations.
\newblock {\em Numer. Math.}, 104\penalty0 (1):\penalty0 103--127, 2006.

\bibitem{S06b}
J.~Szeftel.
\newblock Absorbing boundary conditions for nonlinear scalar partial
  differential equations.
\newblock {\em Comput. Methods Appl. Mech. Engrg.}, 195\penalty0
  (29-32):\penalty0 3760--3775, 2006.

\bibitem{H99}
T.~Hagstrom.
\newblock Radiation boundary conditions for the numerical simulation of waves.
\newblock In {\em Acta numerica, 1999}, volume~8, pages 47--106. Cambridge
  Univ. Press, Cambridge, 1999.

\bibitem{deHoopTD}
A.~de~Hoop, P.~van~den Berg, and R.~Remis.
\newblock Absorbing boundary conditions and perfectly matched layers - an
  analytic time-domain performance analysis.
\newblock {\em IEEE Trans. on Magnetics}, 38:\penalty0 657--660, 2002.

\bibitem{DJ06}
J.~Diaz and P.~Joly.
\newblock A time-domain analysis of {P}{M}{L} models in acoustics.
\newblock {\em Comput. Meth. Appl. Mech. Engrg.}, 195:\penalty0 3820--3853,
  2006.

\bibitem{KC03}
C.A. Kennedy and M.H. Carpenter.
\newblock Additive {Runge-Kutta} schemes for convection-diffusion-reaction
  equations.
\newblock {\em Appl. Numer. Math.}, 44:\penalty0 139--181, 2003.

\bibitem{CM98}
F.~Collino and P.B. Monk.
\newblock {Optimizing the perfectly matched layer}.
\newblock {\em {Comput. Methods Appl. Mech. Engrg.}}, 164:\penalty0 157 -- 171,
  1998.

\end{thebibliography}

\end{document}